\title{Independent sets in tensor graph powers}
\author{{Noga Alon\thanks{ Schools of Mathematics and Computer Science,
Raymond and Beverly Sackler Faculty of Exact Sciences, Tel Aviv
University, Tel Aviv, 69978, Israel. Email: nogaa@tau.ac.il.
Research supported in part by a USA-Israeli BSF grant, by the
Israel Science Foundation and by the Hermann Minkowski Minerva
Center for Geometry at Tel Aviv University.}} \quad {Eyal Lubetzky
\thanks{School of Computer Science, Raymond and Beverly
Sackler Faculty of Exact Sciences, Tel Aviv University, Tel Aviv,
69978, Israel. Email: lubetzky@tau.ac.il. Research partially
supported by a Charles Clore Foundation Fellowship.}}}

\documentclass [11pt]{article}
\usepackage[]{epsf,epsfig,amsmath,amssymb,amsfonts,latexsym,amsthm}

\newtheorem{theorem}{Theorem}[section]
\newtheorem{lemma}[theorem]{Lemma}
\newtheorem{claim}[theorem]{Claim}

\newtheorem{conjecture}[theorem]{Conjecture}
\newtheorem{corollary}[theorem]{Corollary}
\newtheorem{remark}[theorem]{Remark}
\newtheorem{question}[theorem]{Question}
\newtheorem*{question'}{Question \ref{A_is_a_conj}'}
\newtheorem{proposition}[theorem]{Proposition}
\newtheorem{observation}[theorem]{Observation}

\renewcommand{\epsilon}{\varepsilon}

\newcommand{\ignore}[1]{}
\begin{document}
\maketitle

\begin{abstract}
The tensor product of two graphs, $G$ and $H$, has a vertex set
$V(G)\times V(H)$ and an edge between $(u,v)$ and $(u',v')$ iff both
$u u' \in E(G)$ and $v v' \in E(H)$. Let $A(G)$ denote the limit of
the independence ratios of tensor powers of $G$, $\lim
\alpha(G^n)/|V(G^n)|$. This parameter was introduced in \cite{BNR},
where it was shown that $A(G)$ is lower bounded by the vertex
expansion ratio of independent sets of $G$. In this note we study
the relation between these parameters further, and ask whether they
are in fact equal. We present several families of graphs where
equality holds, and discuss the effect the above question has on
various open problems related to tensor graph products.
\end{abstract}

\section{Introduction}
The \textit{tensor} product (also dubbed as categorical or weak
product) of two graphs, $G \times H$, is the graph whose vertex set
is $V(G) \times V(H)$, where two vertices $(u,v)$,$(u',v')$ are
adjacent iff both $u u' \in E(G)$ and $v v' \in E(H)$, i.e., the
vertices are adjacent in each of their coordinates. Clearly, this
product is associative and commutative, thus $G^n$ is well defined
to be the tensor product of $n$ copies of $G$.

The tensor product has attracted a considerable amount of attention
ever since Hedetniemi conjectured in 1966 (\cite{Hedet}) that
$\chi(G \times H) = \min\{\chi(G),\chi(H)\}$ (where $\chi(G)$
denotes the chromatic number of $G$), a problem which remains open
(see \cite{ZhuHedet} for an extensive survey of this problem). For
further work on colorings of tensor products of graphs, see
\cite{ADFS}, \cite{GreenwellLovasz}, \cite{LaroseTardif},
\cite{Tardif}, \cite{TardifChi}, \cite{Zhu}, \cite{ZhuStarChi}.

It is easy to verify that Hedetniemi's conjecture is true when there
is a homomorphism from $G$ to $H$, and in particular when $G=H$, by
examining a copy of $G$ in $G\times H$, and it follows that
$\chi(G^n)=\chi(G)$ for every integer $n$. Furthermore, a similar
argument shows that $\omega(G \times H)$ (the clique number
of $G \times H$) equals
$\min\{\omega(G),\omega(H)\}$ for every two graphs $G$ and $H$, and
in particular, $\omega(G^n) = \omega(G)$ for every integer $n$.
However, the behavior of the independence ratios of the graphs $G^n$ is
far more interesting. Let $i(G)=\alpha(G)/|V(G)|$ denote the
independence ratio of $G$. Notice that for every two graphs $G$ and $H$,
if $I$ is an independent set of $G$,
then the cartesian product $I\times V(H)$ is independent in $G
\times H$, hence every two graphs $G$ and $H$ satisfy:
\begin{equation}\label{i_G_H_geq_max}
  i(G \times H) \geq \max\{i(G),i(H)\}~.
\end{equation}
Therefore, the series $i(G^n)$ is monotone non-decreasing and bounded,
hence its limit exists; we denote this limit, introduced
in \cite{BNR}, where it is called the Ultimate Categorical Independence Ratio of
$G$, by $A(G)$. In contrast to the clique numbers and chromatic
numbers, $A(G)$ may indeed exceed its value at the first power of
$G$, $i(G)$. The authors of \cite{BNR} proved the following simple
lower bound for $A(G)$: if $I$ is an independent set of $G$, then
$A(G) \geq \frac{|I|}{|I|+|N(I)|}$, where $N(I)$ denotes the vertex
neighborhood of $I$. We thus have the following lower bound on
$A(G)$: $A(G) \geq a(G)$, where
$$ a(G) = \max_{\text{$I$ ind. set}} \frac{|I|}{|I|+|N(I)|} ~.$$
It easy to see that $a(G)$ resembles $i(G)$ in the sense that $a(G
\times H) \geq \max\{a(G),a(H)\}$ (to see this, consider the cartesian
product $I \times V(H)$, where $I$ is an independent set of $G$
which attains the ratio $a(G)$). However, as opposed to $i(G)$, it
is not clear if there are any graphs $G,H$ such that $a(G \times H)
> \max\{a(G),a(H)\}$ and yet $a(G),a(H) \leq \frac{1}{2}$. This is
further discussed later.

It is not difficult to see that if $A(G)> \frac{1}{2}$ then
$A(G)=1$, thus $A(G) \in (0,\frac{1}{2}]\cup\{1\}$, as proved in
\cite{BNR} (for the sake of completeness, we will provide short
proofs for this fact and for the fact that $A(G)\geq a(G)$ in
Section \ref{sec::A-equals-a^*}). Hence, we introduce the
following variant of $a(G)$:
$$a^*(G)=\left\{\begin{array}{cl}
a(G) & \mbox{if } a(G) \leq \frac{1}{2} \\
1 & \mbox{if } a(G) > \frac{1}{2}\\
\end{array} \right. ~, $$
and obtain that $A(G) \geq a^*(G)$ for every graph $G$. The
following question seems crucial to the understanding of the
behavior of independence ratios in tensor graph powers: {\question
\label{A_is_a_conj} Does every graph $G$ satisfy $A(G) = a^*(G)$?}

In other words, are non-expanding independent sets of $G$ the only
reason for an increase in the independence ratio of larger powers?
If so, this would immediately settle several open problems related
to $A(G)$ and to fractional colorings of tensor graph products.
Otherwise, an example of a graph $G$ satisfying $A(G) > a^*(G)$
would demonstrate a thus-far unknown way to increase $A(G)$. While
it may seem unreasonable that the complicated parameter $A(G)$
translates into a relatively easy property of $G$, so far the
intermediate results on several conjectures regarding $A(G)$ are
consistent with the consequences of an equality between $A(G)$ and
$a^*(G)$.

As we show later, Question \ref{A_is_a_conj} has the following
simple equivalent form:
\begin{question'}
Does every graph $G$ satisfy $a^*(G^2) = a^*(G)$?
\end{question'}
Conversely, is there a graph $G$ which satisfies the following
two properties:
\begin{enumerate}
\item Every independent set $I$ of $G$ has at least $|I|$ neighbors
(or equivalently, $a(G) \leq \frac{1}{2}$).
\item There exists an independent set $J$ of $G^2$ whose
vertex-expansion ratio, $\frac{|N(J)|}{|J|}$, is strictly
\textit{smaller} than $\frac{|N(I)|}{|I|}$ for every independent set
$I$ of $G$.
\end{enumerate}

In this note, we study the relation between $A(G)$ and $a^*(G)$,
show families of graphs where equality holds, and discuss the
effects of Question \ref{A_is_a_conj} on several conjectures
regarding $A(G)$ and fractional colorings of tensor graph products.
The rest of the paper is organized as follows:

In Section \ref{sec::A-equals-a^*} we present several families of
graphs where equality holds between $A(G)$ and $a^*(G)$. First, we
extend some of the ideas of \cite{BNR} and obtain a characterization
of all graphs $G$ which satisfy the property $A(G)=1$, showing that
for these graphs $a^*(G)$ and $A(G)$ coincide. In the process, we
obtain a polynomial time algorithm for determining whether a graph
$G$ satisfies $A(G)=1$. We conclude the section by observing that
$A(G)=a(G)$ whenever $G$ is vertex transitive, and when
it is the disjoint union of certain vertex transitive graphs.

Section \ref{sec::i-g-times-h} discusses the parameters $i(G)$ and
$a(G)$ when $G$ is a tensor product of two graphs, $G_1$ and $G_2$.
Taking $G_1 = G_2$, we show the equivalence between Questions
\ref{A_is_a_conj} and \ref{A_is_a_conj}'. Next, when $G_1$ and $G_2$
are both vertex transitive, the relation between $i(G)$ and $a(G)$
is related to a fractional version of Hedetniemi's conjecture,
raised by Zhu in \cite{Zhu}. We show that for every two graphs $G$
and $H$, $A(G+H)=A(G \times H)$, where $G+H$ is the disjoint union
of $G$ and $H$. This property links the above problems, along with
Question \ref{A_is_a_conj}, to the problem of determining $A(G+H)$,
raised in \cite{BNR} (where it is conjectured to be equal to
$\max\{A(G),A(H)\}$). Namely, the equality $A(G+H)=A(G \times H)$
implies that if $A(H) = a^*(H)$ for $H=G_1 + G_2$, then:
$$i(G_1 \times G_2) \leq a^*(G_1 + G_2) = \max\{a^*(G_1),a^*(G_2)\}~.$$
This raises the following question, which is a weaker form of
Question \ref{A_is_a_conj}: {\question \label{i_G_H_leq_max_a} Does
the inequality $i(G \times H) \leq \max\{a^*(G),a^*(H)\}$ hold for
every two graphs $G$ and $H$?}

We proceed to demonstrate that several families of graphs satisfy
this inequality, and in the process, obtain several additional
families of graphs $G$ which satisfy $A(G)=a(G)=a^*(G)$.

Section \ref{sec::concluding-remarks} is devoted to concluding
remarks and open problems. We list several additional interesting
questions which are related to $a(G)$, as well as summarize the main
problems which were discussed in the previous sections. Among the
new mentioned problems are those of determining or estimating the value of $A(G)$
for the random graph models $\mathcal{G}_{n,d}$,
$\mathcal{G}_{n,\frac{1}{2}}$ and for the random graph process.

\section{Equality between $A(G)$ and $a^*(G)$}\label{sec::A-equals-a^*}
\subsection{Graphs $G$ which satisfy $A(G)=1$}
In this section we prove a characterization of graphs $G$ satisfying
$A(G)=1$, showing that this is equivalent to the non-existence of a
fractional perfect matching in $G$. A \textit{fractional matching}
in a graph $G=(V,E)$ is a function $f:E\rightarrow \mathbb{R}^+$
such that for every $v \in V$, $\sum_{v \in e}f(e) \leq 1$ (a
matching is the special case of restricting the values of $f$ to
$\{0,1\}$). The value of the fractional matching is defined as $f(E)
= \sum_{e \in E}f(e) ~(\leq \frac{|V|}{2})$. A \textit{fractional
perfect matching} is a fractional matching which achieves this
maximum: $f(E) = \frac{|V|}{2}$.
\begin{theorem}\label{A_of_G_eq_1} For every graph $G$,
$A(G)=1$ iff $a^*(G)=1$ iff $G$
does not contain a fractional perfect matching. \end{theorem}

The proof of Theorem \ref{A_of_G_eq_1} relies on the results of
\cite{BNR} mentioned in the introduction; we recall these results
and provide short proofs for them.
{\claim[\cite{BNR}]\label{A_of_G_geq_a_of_G} For every graph $G$,
$A(G) \geq a(G)$.} {\proof Let $I$ be an independent set which
attains the maximum of $a(G)$. Clearly, for every $k \in
\mathbb{N}$, all vertices in $G^k$, which contain a member of $I
\cup N(I)$ in one of their coordinates, and in addition, whose first
coordinate out of $I \cup N(I)$ belongs to $I$, form an independent
set. As $k$ tends to infinity, almost every vertex has a member of
$I \cup N(I)$ in at least one of its coordinates, and the second
restriction implies that the fractional size of the set above tends
to $\frac{|I|}{|I|+|N(I)|} = a(G)$. \qed }

{\claim [\cite{BNR}]\label{A_of_G_1_2} If $A(G) > \frac{1}{2}$ then
$A(G)=1$.} \begin{proof} Assume, without loss of generality, that
$i(G)>\frac{1}{2}$, and let $I$ be a maximum independent set of $G$.
For every power $k$, the set of all vertices of $G^k$, in which
strictly more than $\frac{k}{2}$ of the coordinates belong to $I$,
is independent. Clearly, since $\frac{|I|}{|G|} > \frac{1}{2}$, the
size of this set tends to $|V(G)|^k$ as $k$ tends to infinity (as the
probability of more Heads than Tails in a sufficiently long sequence
of tosses of a coin biased towards Heads is nearly 1), hence $A(G) =
1$.\end{proof}

\begin{proof}[Proof of Theorem \ref{A_of_G_eq_1}]
By Claims \ref{A_of_G_geq_a_of_G} and \ref{A_of_G_1_2}, if $a(G)
> \frac{1}{2}$ (or equivalently, $a^*(G)=1$) then $A(G) = 1$.
Conversely, assuming that $a^*(G) = a(G) \leq \frac{1}{2}$, we must
show that $A(G)<1$. This will follow from the following simple
lemma, proved by Tutte in 1953 (cf., e.g., \cite{MatchingTheory} p. 216):

\begin{lemma}\label{tutte-lemma}
For a given set $S \subset V(G)$, let $N(S)$ denote that set
of all vertices of $G$ which have a neighbor in $S$; then every set
$S \subset V(G)$ satisfies $|N(S)| \ge |S|$ iff every independent
set $I \subset V(G)$ satisfies $|N(I)| \ge |I|$.
\end{lemma} { \proof [Proof of
lemma] One direction is obvious; for the other direction, take a
subset $S$ with $|N(S)| < |S|$. Define $S'$ to be $\{ v \in S ~|~
N(v) \cap S \neq \emptyset \}$, and examine $I = S \setminus S'$.
Since $S' \subset N(S)$ and $|N(S)| < |S|$, $I$ is
nonempty, and is obviously independent. Therefore $|N(I)| \geq |I|$,
however $|N(I)| \leq |N(S)| - |S'| < |S| - |S'| = |I|$, yielding a
contradiction. \qed}

Returning to the proof of the theorem, observe that by our
assumption that $a(G) \leq \frac{1}{2}$ and the lemma, Hall's
criterion for a perfect matching applies to the bipartite graph $G
\times K_2$ (where $K_2$ is the complete graph on two vertices).
Therefore, $G$ contains a factor $H \subset G$ of vertex disjoint
cycles and edges (to see this, as long as the matching is nonempty,
repeatedly traverse it until closing a cycle and omit these edges).
Since removing edges from $G$ may only increase $A(G)$, it is enough
to show that $A(H)<1$.

We claim that the subgraph $H$ satisfies $A(H) \leq \frac{1}{2}$. To
see this, argue as follows: direct $H$ according to its cycles and
edges (arbitrarily choosing clockwise or counter-clockwise
orientations), and examine the mapping from each vertex to the
following vertex in its cycle. This mapping is an invertible
function $f:V \rightarrow V$, such that for all $v \in V$, $v f(v)
\in E(H)$. Now let $I$ be an independent set of $H^k$. Pick a random
vertex $\underline{u} \in V(H^k)$, uniformly over all the vertices,
and consider the pair $\{\underline{u},\underline{v}\}$, where
$\underline{v}=f(\underline{u})$ is the result of applying $f$ on
each coordinate of $\underline{u}$. Obviously $\underline{v}$ is
uniformly distributed over $H^k$ as-well, thus:
 $$\mathbb{E}\left|I \cap \{\underline{u},\underline{v}\}\right|
 \geq \frac{2}{|H^k|}|I|~.$$
 Choosing a vertex $\underline{u}$ for which $|I \cap \{\underline{u},\underline{v}\}|$
 is at least its expected value, and recalling that $\underline{u}$ and $\underline{v}$
 are adjacent in $H^k$, we get:
  $$ \frac{2}{|H^k|}|I| \leq |I \cap
\{\underline{u},\underline{v}\}|) \leq 1  ~.$$ Hence, $i(H^k) \leq
\frac{1}{2}$, and thus $A(H) \leq \frac{1}{2}$.

An immediate corollary from the above proof that $A(G)=1$ iff $a^*(G)=1$
is the equivalence between the property $A(G) \leq
\frac{1}{2}$ and the existence of a fractional perfect matching in
the graph $G$. It is well known (see for instance
\cite{MatchingTheory}) that for every graph $G$, the maximal
fractional matching of $G$ can be achieved using only the weights
$\{0,\frac{1}{2},1\}$. Therefore, a fractional perfect matching is
precisely a factor $H \subset G$, comprised of vertex disjoint
cycles and edges, and we obtain another format for the condition
$a(G) \leq \frac{1}{2}$: $A(G) \leq
\frac{1}{2}$ iff $G$ has a fractional perfect matching; otherwise,
$A(G) = 1$.

Notice that a fractional perfect matching $f$ of $G$ immediately
induces a fractional perfect matching on $G^k$ for every $k$ (assign
an edge of $G^k$ a weight equaling the product of the weights of
each of the edges in the corresponding coordinates). As it is
easy to see that a fractional perfect matching implies that $i(G) \leq
\frac{1}{2}$, this provides an alternative proof that if $a(G) \leq
\frac{1}{2}$ then $A(G) \leq \frac{1}{2}$.
\end{proof}

Since Lemma \ref{tutte-lemma} also provides us with a polynomial algorithm for
determining whether $a(G)>\frac{1}{2}$ (determine whether Hall's
criterion applies to $G \times K_2$, using network flows), we obtain
the following corollary: {\corollary Given an input graph $G$,
determining whether $A(G)=1$ or $A(G)\leq \frac{1}{2}$ can be done
in polynomial time.}

\subsection{Vertex transitive graphs}
The observation that $A(G)=a(G)$ whenever $G$ is vertex transitive
(notice that $A(G)\leq \frac{1}{2}$ for every nontrivial regular graph $G$) is
a direct corollary of the following result of \cite{ADFS} (the proof
of this fact is by covering $G^k$ uniformly by copies of $G$):
{\proposition [\cite{ADFS}]\label{A_vertex_transitive} If $G$ is
vertex transitive, then $A(G)=i(G)$. }

Clearly, for every graph $G$, $i(G) \leq a(G)$. Hence, for every
vertex transitive graph $G$ the following holds:
$$ A(G) = i(G) \leq a(G) \leq A(G) ~,$$
proving the following corollary:
{\observation\label{a_of_G_vertex_transitive} For every vertex
transitive graph $G$, $A(G)=a^*(G)=a(G)$. }

We conclude this section by mentioning several families of vertex
transitive graphs $G$ and $H$ whose disjoint union $G+H$ satisfies
$A(G+H)=a(G+H)=\max\{A(G),A(H)\}$. These examples satisfy both the
property of Question \ref{A_is_a_conj} and the disjoint union
conjecture of \cite{BNR}.

The next two claims follow from the results of Section
\ref{sec::i-g-times-h}, as we later show. For the first claim,
recall that a circular complete graph (defined in \cite{Zhu}),
$K_{n/d}$, where $n \geq 2d$, has a vertex set $\{0,\ldots,n-1\}$
and an edge between $i,j$ whenever $d \leq |i-j| \leq n-d$. A Kneser
graph, $KN_{n,k}$, where $k \leq n$, has $\binom{n}{k}$ vertices
corresponding to $k$-element subsets of $\{1,\ldots,n\}$, and two
vertices are adjacent iff their corresponding subsets are disjoint.

{\claim\label{disjoint-union-equality-1} Let $G$ and $H$ be two
vertex transitive graphs, where $H$ is one of the following: a
Kneser graph, a circular complete graph, a cycle or a complete
bipartite graph. Then $G+H$ satisfies
$A(G+H)=a(G+H)=\max\{A(G),A(H)\}$. }

{\claim\label{disjoint-union-equality-2} Let $G$ and $H$ be two
vertex transitive graphs satisfying $\chi(G) = \omega(G) \leq
\omega(H)$. Then $A(G+H)=a(G+H)=\max\{A(G),A(H)\}$. }

\section{The tensor product of two graphs}
\label{sec::i-g-times-h}
\subsection{The expansion properties of $G^2$}
Question \ref{A_is_a_conj}, which discusses the relation between the
expansion of independent sets of $G$, and the limit of independence
ratios of tensor powers of $G$, can be translated into a seemingly
simpler question (stated as Question \ref{A_is_a_conj}') comparing
the vertex expansions of a graph and its square: can the minimal
expansion ratio $|N(I)|/|I|$ of independent sets $I$ decrease in the
second power of $G$?

To see the equivalence between Questions \ref{A_is_a_conj} and
\ref{A_is_a_conj}', argue as follows: assuming the answer to
Question \ref{A_is_a_conj} is positive, every graph $G$ satisfies:
$$a^*(G)=A(G) = A(G^2) = a^*(G^2)~,$$ and hence
$a^*(G^2)=a^*(G)$ (recall that every graph $H$ satisfies $a(H^2)
\geq a(H)$). Conversely, suppose that there exists a graph $G$ such
that $A(G) > a^*(G)$. By the simple fact that every
graph $H$ satisfies $i(H) \leq a^*(H)$ we conclude that
there exists an integer $k$ such that
$a^*(G^{2^k}) \geq i(G^{2^k}) > a^*(G)$, and therefore there
exists some integer $\ell \leq k$ for which $a(G^{2^\ell})
> a(G^{2^{\ell-1}})$.

\subsection{The relation between the tensor product and disjoint unions}
In this section we prove the following theorem, which links between
the quantities $i(G_1 \times G_2)$, $a(G_1 \times G_2)$, $\chi_f(G_1
\times G_2)$ and $A(G_1+G_2)$, where $\chi_f(G)$ denotes the
fractional chromatic number of $G$:
\begin{theorem}\label{equiv-statements-vtransitive}For every two vertex
transitive graphs $G_1$ and $G_2$, the following statements are equivalent:
\begin{eqnarray}
\label{i-g1-g2-leq-max-a}
i(G_1 \times G_2) & \leq &
\max\{a^*(G_1),a^*(G_2)\} \\
\label{a-g1-g2-leq-max-a}a^*(G_1 \times G_2) & \leq &
\max\{a^*(G_1),a^*(G_2)\} \\
\label{chi-eq}\chi_f(G_1 \times G_2) & =&
\min\{\chi_f(G_1),\chi_f(G_2)\} \\
\label{A-union-eq} A(G_1 + G_2) &=& \max\{A(G_1),A(G_2)\}
\end{eqnarray}
\end{theorem}
\begin{proof} The proof of Theorem \ref{equiv-statements-vtransitive} relies
on the following proposition: {\proposition
\label{disjoint_union_prod} For every two graphs $G$ and $H$, $A(G +
H) = A(G \times H)$.}

\noindent We note that this generalizes a result of \cite{BNR}, which
states that $A(G+H)$ is at least $\max\{A(G),A(H)\}$. Indeed, that result
immediately follows from the fact that $A(G \times H)$ is always at
least the maximum of $A(G)$ and $A(H)$ (by \eqref{i_G_H_geq_max}).

\begin{proof}[Proof of Proposition \ref{disjoint_union_prod}]
Examine $(G+H)^n$, and observe that a vertex whose $i$-th coordinate
is taken from $G$ is disconnected from all vertices whose $i$-th
coordinate is taken from $H$. Hence, we can break down the $n$-th
power of the disjoint union $G+H$ to $2^n$ disjoint graphs, and
obtain:
\begin{equation}\label{disjoint_break_down}
\alpha\left( (G+H)^n \right) = \sum_{k=0}^n \binom{n}{k}
\alpha\left( G^k H^{n-k} \right)~.
\end{equation}
To prove that $A(G+H) \geq A(G \times H)$, fix $\epsilon
> 0$, and let $N$ denote a sufficiently large integer such that
$i\left( (G \times H)^N\right) \geq (1-\epsilon)A(G \times H)$. The
following is true for every $n
> 2N$ and $N \leq k \leq n-N$:
$$ i(G^k H^{n-k}) = i\left( (G \times H)^N G^{k-N} H^{n-k-N}
\right) \geq i\left((G\times H)^N\right) \geq (1-\epsilon)A(G\times
H)~,$$ where the first inequality is by \eqref{i_G_H_geq_max}. Using
this inequality together with \eqref{disjoint_break_down} yields:
\begin{eqnarray}
i\left( (G + H)^n \right) & \geq & \frac{1}{|(G + H)^n|}
\sum_{k=N}^{n-N} \binom{n}{k} \alpha(G^k H^{n-k}) \geq \nonumber \\
& \geq & \frac{\sum_{k=N}^{n-N} \binom{n}{k} |G|^k
|H|^{n-k}}{(|G|+|H|)^n} (1-\epsilon) A(G\times H)
\mathop{\longrightarrow}_{n \rightarrow \infty} (1-\epsilon)A(G
\times H) ~.\nonumber
\end{eqnarray} Therefore $A(G + H) \geq (1-\epsilon)A(G \times H)$ for any
$\epsilon > 0$, as required.

It remains to show that $A(G+H) \leq A(G \times H)$. First
observe that \eqref{i_G_H_geq_max} gives the following
relation:
\begin{equation}\label{monotonicity_i_of_G} \forall ~k,l \geq 1
~,~~ i(G^k H^l) \leq i(G^k H^l \times G^l H^k) = i(G^{k+l} H^{k+l})
\leq A(G \times H)~.
\end{equation} Using \eqref{disjoint_break_down} again, we obtain:
\begin{eqnarray}
 i \left( (G + H)^n \right) &=& \sum_{k=0}^n \binom{n}{k}
\frac{\alpha(G^k H^{n-k})}{(|G|+|H|)^n} = \sum_{k=0}^n \binom{n}{k}
i(G^k H^{n-k}) \cdot \frac{|G|^k |H|^{n-k}}{(|G|+|H|)^n}
\nonumber \\
&\leq& \frac{|G|^n}{(|G|+|H|)^n} i(G^n) + \frac{|H|^n}{(|G|+|H|)^n}
i(H^n) + \left( 1 - \frac{|G|^n+|H|^n}{(|G|+|H|)^n} \right)A(G
\times H) \nonumber \\
&\leq& \frac{|G|^n}{(|G|+|H|)^n} A(G) + \frac{|H|^n}{(|G|+|H|)^n}
A(H) + \left( 1 - \frac{|G|^n+|H|^n}{(|G|+|H|)^n} \right)A(G \times
H) \nonumber \\
&\longrightarrow& A(G \times H) ~,\nonumber
\end{eqnarray}
where the first inequality is by \eqref{monotonicity_i_of_G}, and
the second is by definition of $A(G)$.\end{proof}

Equipped with the last proposition, we can now prove that Question
\ref{i_G_H_leq_max_a} is indeed a weaker form of Question
\ref{A_is_a_conj}, namely that if $A(G)=a^*(G)$ for every $G$, then
$i(G_1 \times G_2) \leq \max\{a^*(G_1),a^*(G_2)\}$ for every two
graphs $G_1,G_2$. Indeed, if $A(G_1+G_2)=a^*(G_1+G_2)$ then
inequality \eqref{i-g1-g2-leq-max-a} holds, as well as the stronger
inequality \eqref{a-g1-g2-leq-max-a}:
$$
a^*(G_1 \times G_2) \leq A(G_1 \times G_2) = A(G_1 + G_2) = a^*(G_1 +
G_2) = \max\{a^*(G_1),a^*(G_2)\}~,$$ as required.

Having shown that a positive answer to Question \ref{A_is_a_conj} implies
inequality \eqref{a-g1-g2-leq-max-a} (and hence inequality \eqref{i-g1-g2-leq-max-a}
as well), we show the implications of inequality \eqref{i-g1-g2-leq-max-a} when
the two graphs are vertex transitive.

Recall that for every two graphs $G_1$ and $G_2$, $i(G_1 \times G_2)
\geq \max\{i(G_1),i(G_2)\}$, and consider the case when $G_1,G_2$
are both vertex transitive and have edges. In this case,
$i(G_i)=a(G_i)=a^*(G_i)$ ($i=1,2$), hence inequalities
\eqref{i-g1-g2-leq-max-a} and \eqref{a-g1-g2-leq-max-a} are
equivalent, and are both translated into the form $i(G_1 \times G_2)
= \max\{i(G_1),i(G_2)\}$. Next, recall that for every vertex
transitive graph $G$, $i(G) = 1/\chi_f(G)$. Hence, inequality
\eqref{i-g1-g2-leq-max-a} (corresponding to Question
\ref{i_G_H_leq_max_a}), when restricted to vertex transitive graphs,
coincides with \eqref{chi-eq}. Furthermore, by Observation
\ref{a_of_G_vertex_transitive} and Proposition
\ref{disjoint_union_prod}, for vertex transitive $G_1$ and $G_2$ we
have:
$$i(G_1 \times G_2) = A(G_1 \times G_2) = A(G_1 + G_2) \geq
\max\{A(G_1),A(G_2)\} = \max\{i(G_1),i(G_2)\}~,$$
hence in this case \eqref{chi-eq} also coincides with \eqref{A-union-eq}.
Thus, all four statements are equivalent for vertex transitive graphs. \end{proof}

By the last theorem, the following two conjectures, raised in
$\cite{BNR}$ and $\cite{Zhu}$, coincide for vertex transitive
graphs:
{\conjecture[\cite{BNR}]\label{c33} For every two graphs $G$ and $H$,
$A(G+H)=\max\{A(G),A(H)\}$.}
{\conjecture[\cite{Zhu}]\label{chi_f_conj}For every two graphs $G$
and $H$, $\chi_f(G \times H)=\min\{\chi_f(G),\chi_f(H)\}$.}

The study of Conjecture \ref{chi_f_conj} is somewhat related to the
famous and long studied Hedetniemi conjecture (stating that $\chi(G
\times H) = \min\{\chi(G),\chi(H)\}$), as for every two graphs $G$
and $H$, $\omega(G \times H) = \min\{\omega(G),\omega(H)\}$, and
furthermore $\omega(G) \leq \chi_f(G) \leq \chi(G)$.

It is easy to see that the inequality $\chi_f(G\times H) \leq
\min\{\chi_f(G),\chi_f(H)\}$ is always true. It is shown in
\cite{Tardif} that Conjecture \ref{chi_f_conj} is not far from being
true, by proving that for every graphs $G$ and $H$,
$\chi_f(G \times H) \geq \frac{1}{4}
\min\{\chi_f(G),\chi_f(H)\}$.

So far, Conjecture \ref{chi_f_conj} was verified (in \cite{Zhu}) for
the cases in which one of the two graphs is either a Kneser graph or
a circular-complete graph. This implies the cases of $H$ belonging
to these two families of graphs in Claim
\ref{disjoint-union-equality-1}. Claim
\ref{disjoint-union-equality-2} is derived from the the following
remark, which provides another family of graphs for which Conjecture
\ref{chi_f_conj} holds.

 {\remark Let $G$ and $H$ be graphs such that
$\chi(G)=\omega(G) \leq \omega(H)$. It follows that $\omega(G \times
H) = \omega(G) = \chi(G \times H)$, and thus $\chi_f(G \times H)
=\min\{\chi_f(G),\chi_f(H)\}$, and $\chi(G \times H)
=\min\{\chi(G),\chi(H)\}$. In particular, this is true when $G$ and
$H$ are perfect graphs. }

\subsection{Graphs satisfying the property of Question
\ref{i_G_H_leq_max_a}} In this subsection we note that several
families of graphs satisfy inequality \eqref{i-g1-g2-leq-max-a} (and
the property of Question \ref{i_G_H_leq_max_a}). This appears in the
following propositions:

{\proposition \label{i_of_G_times_C_l} For every graph $G$ and
integer $\ell$, $i(G \times C_\ell) \leq \max\{ a(G), a(C_\ell\}$
(and hence, $i(G \times C_\ell) \leq \max\{a^*(G),a^*(C_\ell)\}$).
This result can be extended to $G \times H$, where $H$ is a disjoint
union of cycles.}
\begin{proof} We need the following lemma:
\begin{lemma}\label{full-copy-lemma} Let $G$ and $H$ be two graphs which satisfy at least
one of the following conditions: \\
\indent 1. $a(G) \geq \frac{1}{2}$, and every $S \subsetneqq V(H)$
satisfies
$|N(S)|>|S|$. \\
\indent 2. $a(G) > \frac{1}{2}$ and every $S \subset V(H)$ satisfies
$|N(S)| \geq |S|$. \\
Then every maximum independent set $I \subset V(G \times H)$
contains at least one ``full copy'' of $H$, i.e., for each such $I$
there is a vertex $v \in V(G)$, such that $\{(v,w) : w \in H\}
\subset I$.\end{lemma}
\begin{proof}[Proof of lemma] We begin with
the case $a(G) \geq \frac{1}{2}$ and $|N(S)| > |S|$ for every $S
\subsetneqq V(H)$. Let $J$ be a smallest (with respect to either
size or containment) nonempty independent set in $G$ such that
$\frac{|J|}{|J|+|N(J)|} \geq \frac{1}{2}$. Clearly, $|N(J)| \leq
|J|$. We claim that this inequality proves the existence of a
one-one function $f:N(J) \rightarrow J$, such that $v f(v) \in E(G)$
(that is, there is a matching between $N(J)$ and $J$ which saturates
$N(J)$). To prove this fact, take any set $S \subset N(J)$ and
assume $|N(S) \cap J| < |S|$; it is thus possible to delete $N(S)
\cap J$ from $J$ (and at least $|S|$ vertices from $N(J)$) and since
$|N(S) \cap J| < |S| \leq |N(J)| \leq |J|$ we are left with a
nonempty $J' \subsetneqq J$ satisfying $|N(J')| \leq |J'|$. This
contradicts the minimality of $J$. Now we can apply Hall's Theorem
to match a unique vertex in $J$ for each vertex in $N(J)$.

Assume the lemma is false, and let $I$ be a counterexample. Examine
the intersection of $I$ with a pair of copies of $H$, which are
matched in the matching above between $N(J)$ and $J$. As we assumed
that there are no full $H$ copies in $I$, each set $S$ of vertices
in a copy of $H$ has at least $|S|+1$ neighbors in an adjacent copy
of $H$. Thus, each of the matched pairs of $N(J) \rightarrow J$
contains at most $|H|-1$ vertices of $I$. Define $I'$ as the result
of adding all missing vertices from the $H$ copies of $J$ to $I$,
and removing all existing vertices from the copies of $N(J)$ (all
other vertices remain unchanged). Then $I'$ is independent, and we
obtain a contradiction to the maximality of $I$.

The case of $a(G) > \frac{1}{2}$ and $|N(S)| \geq |S|$ for every $S
\subset V(H)$ is essentially the same. The set $J$ is now the
smallest independent set of $G$ for which $|N(J)| < |J|$, and again,
this implies the existence of a matching from $N(J)$ to $J$, which
saturates $N(J)$. By our assumption on $H$, each pair of copies of
$H$ in the matching contributes at most $|H|$ vertices to a maximum
independent set $I$ of $G \times H$, and by the assumption on $I$,
the unmatched copies of $H$ (recall $|J| > |N(J)|$) are incomplete.
Therefore, we have strictly less than $|H||J|$ vertices of $I$ in
$(N(J) \cup J)\times H$, contradicting the maximality of
$I$.\end{proof}

Returning to the proof of the proposition, let $I$ be a maximum
independent set of $G \times C_\ell$. Remove all vertices, which
belong to full copies of $C_\ell$ in $I$, if there are any, along
with all their neighbors (note that these neighbors are also
complete copies of $C_\ell$, but this time empty ones). These
vertices contribute a ratio of at most $a(G)$, since their copies
form an independent set in $G$. Let $G'$ denote the induced graph of
$G$ on all remaining copies. The set $I'$, defined to be $I \cap (G'
\times C_\ell)$, is a maximum independent set of $G' \times C_\ell$,
because for any member of $I$ we removed, we also removed all of its
neighbors from the graph.

Notice that $C_\ell$ satisfies the expansion property required from
$H$ in Lemma \ref{full-copy-lemma}: for every $k$, every set $S
\subsetneqq V(C_{2k+1})$ satisfies $|N(S)| > |S|$, and every set $S
\subset V(C_{2k})$ satisfies $|N(S)| \geq |S|$. We note that, in
fact, by the method used in the proof of Lemma \ref{tutte-lemma} it
is easy to show that every regular graph $H$ satisfies $|N(S)| \geq
|S|$ for every set $S\subset V(H)$, and if in addition $H$ is
non-bipartite and connected, then every $S \subsetneqq V(H)$
satisfies $|N(S)|>|S|$.

We can therefore apply the lemma on $G' \times C_\ell$. By
definition, there are no full copies of $C_\ell$ in $I'$, hence, by
the lemma, we obtain that $a(G') \leq \frac{1}{2}$ (and even $a(G')
< \frac{1}{2}$ in case $\ell$ is odd). In particular, we can apply
Hall's Theorem and obtain a factor of edges and cycles in $G'$. Each
connected pair of non-full copies has an independence ratio of at
most $i(C_\ell) = a(C_\ell)$ (by a similar argument to the one
stated in the proof of the lemma), and double counting the
contribution of the copies in the cycles we conclude that
$\frac{|I'|}{|G'||C_\ell|} \leq a(C_\ell)$. Therefore $i(G \times
C_\ell)$ is an average between values, each of which is at most
$\max\{a(G),a(C_\ell)\}$, completing the proof. \end{proof}

{\proposition \label{i_of_G_times_K_k} For every graph $G$ and
integer $k$: $i(G \times K_k) \leq \max\{ a(G), a(K_k)\}$ (and
hence, such graphs satisfy the inequality of Question
\ref{i_G_H_leq_max_a}). This result can be extended to $G \times H$,
where $H$ is a disjoint union of complete graphs.}

\begin{proof}
\ignore{ Let us start with the case $a(G) \leq a(K_k) =
\frac{1}{k}$. Let $I$ be a maximum independent set of $G\times K_k$,
and let $S_i$ ($i \in \{1,\ldots,k\}$) denote its intersection with
the $i$-th copy of $G$. It immediately follows that there are no
edges of $G$ between $S_i$ and $S_j$ for all $i \neq j$. Assume that
$i(G \times K_k)
> \frac{1}{k}$; then $|I| = \sum_i |S_i| > |G|$. Now, suppose that some
vertex $v \in V(G)$ belongs to $S_i \cap S_j$ for $i \neq j$. This
implies that this vertex is disconnected from all the sets $S_\ell$,
$\ell\in \{1,\ldots,k\}$, thus we could add it to any group $S_\ell$
which did not include it in the first place. Since $I$ is maximum,
such a $v$ indeed belongs to all groups $S_\ell$. Let $X = \bigcap_i
S_i$, and set $S'_i = S_i \setminus X$. The groups
$S'_1,\ldots,S'_k,X$ are pairwise disjoint and pairwise
disconnected, and in addition, $X$ is independent. Finally:
\begin{equation}\label{k_X_eq} |I| = \sum_i |S_i| = k |X| +
\sum_i |S'_i| > |G| ~.\end{equation} We now examine the neighbor
ratio of $X$ in $G$: all the potential neighbors of $X$ belong to
the set $V(G) \setminus \left(S'_1 \cup \ldots \cup S'_k \cup
X\right) $, therefore $$ a(G) \geq \frac{|X|}{|X|+|N(X)|} \geq
\frac{|X|}{|X| + (|G| - |S'_1| - \ldots - |S'_k| - |X|)} =
\frac{|X|}{|G| - \sum_i|S'_i|} > \frac{1}{k} ~,$$ where the last
inequality results from equation \eqref{k_X_eq}. This contradicts
the assumption that $a(G) \leq \frac{1}{k}$.

We are left with the case $a(G) > \frac{1}{k}$. } Let $I$ denote a
maximum independent set of $G \times K_k$, and examine all copies of
$K_k$ which contain at least two vertices of $I$. Such a copy of
$K_k$ in $G\times K_k$ forces its neighbor copies to be empty (since
two vertices of $K_k$ have the entire graph $K_k$ as their
neighborhood). Therefore, by the maximality of $I$, such copies must
contain \textit{all} vertices of $K_k$. Denote the vertices of $G$
which represent these copies by $S \subset V(G)$; then $S$ is an
independent set of $G$, and the copies represented by $S \cup N(S)$
contribute an independence ratio of at most $a(G)$. Each of the
remaining copies contains at most one vertex, giving an independence
ratio of at most $\frac{1}{k}=a(K_k)$. Therefore, $i(G \times K_k)$
is an average between values which are at most
$\max\{a(G),a(K_k)\}$, and the result follows.
\end{proof}

{\corollary \label{i_of_a(G)=1/2} Let $G$ be a graph satisfying
$a(G) = \frac{1}{2}$; then for every graph $H$ the following
inequality holds: $i(G \times H) \leq \max\{ a(G), a(H)\}$.}

\begin{proof} By Theorem \ref{A_of_G_eq_1} we deduce that $G$ contains a
fractional perfect matching; let $G'$ be a factor of $G$ consisting
of vertex disjoint cycles and edges. Since $a(G') \leq \frac{1}{2}$
as-well, it is enough to show that $i(G' \times H) \leq \max\{
a(G'), a(H)\}$. Indeed, since $G'$ is a disjoint union of the form
$C_{\ell_1} + \ldots + C_{\ell_k} + K_2 + \ldots + K_2$, the result
follows from Proposition \ref{i_of_G_times_C_l} and Proposition
\ref{i_of_G_times_K_k}.
\end{proof}

\section{Concluding remarks and open problems}\label{sec::concluding-remarks}

We have seen that answering Question \ref{A_is_a_conj} is imperative
to the understanding of the behavior of independent sets in tensor
graph powers. While it is relatively simple to show that $A(G)$
equals $a(G)$ whenever $G$ is vertex transitive, proving this
equality for $G=G_1+G_2$, the disjoint union of two vertex
transitive graphs $G_1$ and $G_2$, seems difficult; it is equivalent
to Conjecture \ref{chi_f_conj}, the fractional version of
Hedetniemi's conjecture, for vertex transitive graphs. These two
conjectures are consistent with a positive answer to Question
\ref{A_is_a_conj}, and are in fact direct corollaries in such a
case.

The assertion of Conjecture \ref{c33} for several cases can be deduced
from the spectral bound for $A(G)$ proved in \cite{ADFS}. For a regular
graph $G$ with $n$ vertices and
eigenvalues $\lambda_1 \geq \lambda_2 \geq \cdots \geq \lambda_n$,
denote $\Lambda(G)=\frac{-\lambda_n}{\lambda_1 -\lambda_n}$.
As observed in \cite{ADFS}, the usual known spectral upper bound for the independence
number of a graph implies that for each regular $G$, $A(G) \leq \Lambda(G)$.
It is not difficult to check that for regular $G$ and $H$,
$\Lambda(G \times H)=\max\{\Lambda(G),\Lambda(H)\}$. Therefore,
by Proposition \ref{disjoint_union_prod}, if $G$ and $H$ are regular and satisfy
$\Lambda(G) \leq \Lambda(H)=A(H)$, then the assertion of Conjecture \ref{c33}
holds for $G$ and $H$. Several examples of graphs $H$ satisfying
$\Lambda(H)=A(H)$ are mentioned in \cite{ADFS}.

It is interesting to inspect the expected values of $A(G)$ for
random graph models. First, consider $G^t \sim {\cal G}_{n,t}$, the
random graph process on $n$ vertices after $t$ steps, where there
are $t$ edges chosen uniformly out of all possible edges (for more
information on the random graph process, see \cite{RandomGraphs}).
It is not difficult to show, as mentioned in \cite{Isoperimetric},
that $a(G^t)$ equals the minimal degree of $G^t$, $\delta(G^t)$, as
long as $\delta(G^t)$ is fixed and $|G|$ is sufficiently large. When
considering $A(G)$, the following is a direct corollary of the
fractional perfect matching characterization for $A(G)=1$ (Theorem
\ref{A_is_a_conj}), along with the fact that the property ``$G$
contains a fractional perfect matching" almost surely has the same
hitting time as the property ``$\delta(G) \geq 1$":

{\remark With high probability, the hitting time of the property
$A(G)<1$
 equals the hitting time of $\delta(G)\geq 1$.
 Furthermore, almost every graph process at that time satisfies
 $A(G)=\frac{1}{2}$. }

 {\question Does almost every graph process satisfy
 $A(G)=\frac{1}{\delta(G)+1}$ as long as $\delta(G)$ is fixed? }

Second, the expected value of $a(G)$ for a random regular graph $G \sim
\mathcal{G}_{n,d}$ is easily shown to be $\Theta(\frac{\log d}{d})$, as
the independence ratio of $\mathcal{G}_{n,d}$ is almost surely
between $\frac{\log d}{d}$ and $2\frac{\log d}{d}$ as
$n \rightarrow \infty$ (see \cite{RandomRegularUpperRatio},
\cite{RandomRegularLowerRatio}). As for $A(G)$,
the following is easy to prove, by the spectral upper bound $\Lambda(G)$ mentioned
above,
and by the eigenvalue estimations of \cite{RandomRegularEigenvalues}:
 {\remark Let $G$ denote the random regular graph
$\mathcal{G}_{n,d}$; almost surely: $\Omega(\frac{\log
d}{d}) \leq A(G) \leq O(\frac{1}{\sqrt{d}})$ as $d\rightarrow \infty$.}

{\question Is the expected value of $A(G)$ for the random regular
graph $G \sim \mathcal{G}_{n,d}$ equal to $\Theta(\frac{\log d}{d})$?}

The last approach can be applied to the random graph $G \sim
\mathcal{G}_{n,\frac{1}{2}}$ as well. To see this, consider
a large regular factor (see \cite{ShamirUpfal}), and use the eigenvalue estimations of
\cite{RandomEigenvalues} to obtain that almost surely $\Omega(\frac{\log
n}{n}) \leq A(G) \leq O(\sqrt{\frac{\log n}{n}})$, whereas $a(G)$ is
almost surely $(2+o(1))\frac{\log_2 n}{n}$.

{\question Is the expected value of $A(G)$ for the random graph $G
\sim \mathcal{G}_{n,\frac{1}{2}}$ equal to $\Theta(\frac{\log n}{n})$? }

We conclude with the question of the decidability of $A(G)$.
Clearly, deciding if $a(G) > \beta$ for a given value $\beta$ is in
NP, and we can show that it is in fact NP-complete. It seems
plausible that $A(G)$ can be calculated (though not necessarily by
an efficient algorithm) up to an arbitrary precision: {\question Is
the problem of deciding whether $A(G) > \beta$, for a given graph
$G$ and a given value $\beta$, decidable? }

\end{document}